\documentclass[a4paper,11pt]{article}
\usepackage{latexsym}
\usepackage{amsfonts}
\usepackage{amssymb,amsmath}
\usepackage{amscd}
\usepackage[dvips]{graphicx}
\usepackage{color}
\usepackage{amsthm}
\usepackage{fullpage}

\author{ M. Dziemia\'nczuk\footnote{
Institute of Informatics, University of Gda\'nsk, Poland;  E-mail: \texttt{mdziemianczuk@gmail.com}.}
}
\title{\textsc{Generalization of Fibonomial Coefficients}}

\newtheorem{theoremn}{Theorem}

\newtheorem{proposition}{Proposition}

\newtheorem{corol}{Corollary}

\newenvironment{Proof}{\noindent{\it Proof.}}{\nolinebreak[4]\hfill$\blacksquare$\\\par}

\theoremstyle{definition}
\newtheorem{defn}{Definition}[section]
\newtheorem{exmp}{Example}[section]

\numberwithin{equation}{section}

\newcommand{\layer}[2]{\nobreak{\langle\Phi_{#1}\!\to\!\Phi_{#2}\rangle}}

\newcommand{\fnomial}[2]{\binom{\,#1\,}{#2}_{\!\!F} }
\newcommand{\fnomialF}[3]{\binom{\,#1\,}{#2}_{\!\!#3} }
\newcommand{\fnomialT}[2]{\binom{\,#1\,}{#2}_{\!\mathcal{T}}}

\newcommand{\TseqS}[1]{\mathcal{T}_{#1}}
\newcommand{\TseqPQ}[2]{\mathcal{T}_{#1,#2}}
\newcommand{\Tseq}[0]{\mathcal{T}}

\newcommand{\Nat}{\mathbb{N}}

\begin{document}
%\maketitle

\begin{center}
	{\LARGE\textsc{Generalization of Fibonomial Coefficients}}
	
	\vspace{0.6cm}
	
	{\Large M. Dziemia\'nczuk}\footnote{Institute of Informatics, University of Gda\'nsk, Poland;  E-mail: \texttt{mdziemianczuk@gmail.com}.}
	\vspace{1.4cm}
\end{center}

\begin{abstract}
Following Lucas and then other Fibonacci people Kwa\'sniewski had introduced and had started ten years ago the open investigation of the overall $F$-nomial coefficients which encompass among others Binomial, Gaussian and Fibonomial coefficients with a new unified combinatorial interpretation expressed in terms of cobweb posets' partitions and tilings of discrete hyperboxes. In this paper, we deal with special subfamily of $\Tseq$-nomial coefficients.

	The main aim of this note is to develop the theory of $\Tseq$-nomial coefficients with the help of generating functions. The binomial-like theorem for $\mathcal{T}$-nomials is delivered here and some consequences  of it are drawn. A new combinatorial interpretation of $\Tseq$-nomial coefficients is provided and compared with the Konvalina way of objects' selections from weighted boxes. A brief summary of already known properties of $F$-nomial coefficients is served.

\vspace{0.4cm}

\noindent  This is The Internet Gian-Carlo Rota Polish Seminar article, No \textbf{9}, \textbf{Subject 5}, \textbf{2009-08-08}, \emph{http://ii.uwb.edu.pl/akk/sem/sem\_rota.htm}\\

\end{abstract}

% % % % % % % % % % % % % % % % % % % % % % % % % % % % % % % % % % 
\section{Introduction}

At first let us recall definition of the \emph{$F$-nomial} coefficients and summarize already known properties of these arrays of nonnegative integer numbers. Next, a special family of tileable sequences $\Tseq$ and its corresponding \emph{$\Tseq$-nomial} coefficients is being considered.

\begin{defn}\label{def:fnomial}
Let $F$ stays for a natural numbers' sequence $\{n_F\}_{n\geq 0}$, and $n,k\in\Nat$, such that $n\geq k$. Then \emph{$F$-nomial} coefficient is identified with the symbol
\begin{equation}
	\fnomial{n}{k} = \frac{n_F!}{k_F!(n-k)_F!} = \frac{n_F^{\underline{k}}}{k_F!}
\end{equation}
where $n_F! = n_F (n-1)_F\cdots 1_F$ and $n_F^{\underline k} = n_F (n-1)_F\cdots (n-k+1)_F$ with $0_F! = n_F^{\underline 0} = 1$.
\end{defn}

Since $F\equiv\{n_F\}_{n\geq 1}$ stays for a sequence of natural numbers, i.e. $n_F = n$, the \emph{$F$-nomial} coefficient reduce to ordinary binomial coefficient

\begin{equation*}
	\fnomial{n}{k} =  \frac{n!}{k!(n-k)!} = \fnomialF{n}{k}{}
\end{equation*}
For a sequence $F$ of next Fibonacci numbers $\{\mathcal{F}_n\}_{n\geq 0}$ we obtain Fibonomial coefficient, i.e.
\begin{equation*}
	\fnomial{n}{k} = \frac{\mathcal{F}_n!}{\mathcal{F}_k! \mathcal{F}_{n-k}!} = \fnomialF{n}{k}{Fib}
\end{equation*}
Finally, if an $n$-th element of the sequence $F$ is $n_F = n_q = (q^n - 1)/(q - 1)$ we obtain $q$-binomial (Gaussian) coefficient 
\begin{equation*}
	\fnomial{n}{k} = \frac{n_q!}{k_q!(n-k)_q!} = \fnomialF{n}{k}{q}.
\end{equation*}

\noindent 
Let us review main properties of \emph{$F$-nomial} coefficients, where $F$ denote any natural numbers' valued sequence. According to the definition, we have

\begin{enumerate}

\item \emph{Complementation Rule}
\begin{equation}
	\label{eq:symmetry}
	\fnomial{n}{k} = \fnomial{n}{n-k},
\end{equation}

\item \emph{Iterative Rule}
\begin{equation}
	\fnomial{n}{m}\fnomial{m}{k} = \fnomial{n}{k}\fnomial{n-k}{n-m}
\end{equation}

\item \emph{Multinomial coefficients} 
\begin{equation}
	\fnomial{n}{i_1,i_2,\ldots,i_k} = \fnomial{n}{i_1}\fnomial{n-i_1}{i_2}\cdots\fnomial{n-i_1-\cdots-i_{k-1}}{i_k}
\end{equation}

\item \emph{Inversion formula} 
\begin{equation}
	\label{eq:inv}
	\fnomial{n}{k}^{-1} = \fnomial{n}{k}\sum_{s=1}^{n-k} (-1)^s \sum_{\substack{i_1+i_2+\cdots+i_s=n-k \\ i_1,i_2,\ldots,i_s\geq 1}} \fnomial{n-k}{i_1,i_2,\ldots,i_s}
\end{equation}
while $\binom{n}{n}_F^{-1} = 1$.

\end{enumerate}

Since a special  $\mathcal{A}$-\emph{admissible} family of natural numbers' valued sequences $F$, introduced by Kwa\'sniewski \cite{akk6,akk1,akk2} is taken into account, the \emph{$F$-nomial} coefficients counts blocks of \emph{cobweb poset}'s partitions. This family includes for example Natural numbers, Gaussian and Fibonacci integers.

\vspace{0.2cm}
\textbf{Combinatorial Interpretation I} (``\emph{Partitions of cobweb layer}'')

\vspace{0.2cm}
\noindent
Let $\layer{k}{n}$ be a cobweb poset layer with $m$ levels $\Phi$, where $k=n-m+1$, and $C_{max}(\Pi)$ denote the number of maximal chains of a poset $\Pi$. Suppose that $F$ is an \emph{admissible} sequence $\mathcal{A}$. Then
\begin{equation}
	\layer{n}{k} = \bigcup_{i=1}^{\lambda} \pi_i
	\quad \Leftrightarrow \quad
	\lambda = \fnomial{n}{m}
\end{equation}
while $C_{max}(\pi_i) = C_{max}(\layer{1}{m})$ and $\pi_i\cap\pi_j = \emptyset$ for any $i\neq j$.

\vspace{0.2cm}
Kwa\'sniewski posed also the \emph{cobweb tiling} problem \cite{akk1}, where one asks about family of so-called \emph{tileable} sequences $\mathcal{T}$. For such sequences, the \emph{$F$-nomial} coefficients obtain additional combinatorial interpretation with respect to the general Interpretation I.

\vspace{0.2cm}
\textbf{ Combinatorial Interpretation II} (``\emph{Tilings of hyper $F$-boxes}'')

\vspace{0.2cm}
	Let $V_{k,n}$ denote an $m$-dimensional discrete box $A_1\times A_2\times \cdots \times A_m$ where $|A_s| = (k+s-1)_F$ and $m = n-k+1$ (see: \cite{akk5,md1,akk7}).
Suppose $F$ is a \emph{tileable} sequence $\mathcal{T}$. Then the value of \emph{$\Tseq$-nomial} coefficient $\binom{n}{m}_{\Tseq}$ is equal to the number of $m$-dimensional translates bricks $V_{1,m}$ that form a tiling of $V_{k,n}$.

\vspace{0.2cm}
For further reading about combinatorial interpretations we refer the reader to \cite{akk1,akk6,akk7,md1,md2} and to the references therein.

% % % % % % % % % % % % % % % % % % % % % % % % % % % % % % % % % % 
\section{Tileable sequences and $\Tseq$-nomial coefficients}

The Kwa\'sniewski \emph{upside-down} notation \cite{akk1,akk2}, invented in the spirit of Knuth \cite{knuth} is being here taken for granted. For example $n$-th element of a sequence $F=\{n_F\}_{n\geq 1}$ is $F_n \equiv n_F$, consequently $n_F! = n_F \cdot (n - 1)_F \cdots 1_F$, and $n_F^{\underline{k}} = n_F\cdot(n-1)_F\cdots(n-k+1)_F$ with $n_F^{\underline{0}} = 0_F! = 1$.

\vspace{0.2cm}
In this section we define so-called tileable sequences $\TseqPQ{p}{q}$ for arbitrary $p$ and $q$. In the next part of this paper, the corresponding \emph{$\Tseq$-nomial} coefficients, determined by these sequences are considered. 

\begin{defn}
A natural numbers' valued sequence $\Tseq\equiv\{n_{\Tseq}\}_{n\geq 1}$ constituted by $n$-th coefficient of the generating function $\TseqPQ{p}{q}(x)$ expansion, i.e., $n_{\Tseq} = [x^n]\TseqPQ{p}{q}(x)$, where
\begin{equation}\label{eq:form}
	\TseqPQ{p}{q}(x) = 1_{\Tseq} \cdot \frac{x}{(1 - p\,x)(1 - q\,x)}
\end{equation}
while $p,q\in \mathbb{R}$ and $1_{\Tseq} \in \mathbb{R}$ is called \emph{tileable} and denoted by $\TseqPQ{p}{q}$.
\end{defn}

Take a \emph{tileable} sequence $\Tseq \equiv \TseqPQ{p}{q}$. Without loss of generality we assume that $1_{\Tseq} = 1$. Let $n\in\Nat$ be given. Then for any $m,k\in\Nat$ such that $m+k = n$ an $n$-th element of $\Tseq$ satisfies the following recurrence relation
\begin{equation}
	\label{eq:recurseq}
	n_{\Tseq} = (k+m)_{\Tseq} = p^m k_{\Tseq} + q^k m_{\Tseq}.
\end{equation}
We may generalize the above as follows. Let $\mathbf{b}$ be a composition $\langle b_1,b_2,\ldots,b_k\rangle$ of the number $n$ into $k$ non-zero parts. Then an $n$-th element of the sequence $\Tseq$ satisfies
\begin{equation}
	\label{eq:recurseqmulti}
	n_{\Tseq} = \big( b_1+b_2+\cdots+b_k \big)_{\!\Tseq} = \sum_{i=1}^k  p^{(b_{i+1} + \cdots + b_{k})} q^{(b_1+ \cdots +b_{i-1})}  (b_i)_{\Tseq}.
\end{equation}

\noindent An explicit formula for $n$-th term of $\Tseq$ is given by
\begin{equation}
	n_{\Tseq} =
	\left\{ 
	\begin{array}{ll}
		\frac{q^n - p^n}{q-p}, &  q\neq p,
		\\
	\\
		n\, q^{n-1}, & q = p.
	\end{array}
	\right.
\end{equation}
for $n\geq 1$. In the other hand we have
\begin{equation}
	n_{\Tseq} = \sum_{i=1}^{n} q^{(n-i)} p^{(i-1)}.
\end{equation}
The explicit formula of \emph{$\Tseq$-nomial} coefficients while the sequence $\Tseq = \TseqPQ{p}{q}$ is
\begin{equation}
	\fnomialT{n}{k} = 
	\left\{ 
	\begin{array}{ll}
	\prod_{i=1}^{k} \frac{(p^{n-i+1} - q^{n-i+1})}{(p^i-q^i)}, & p \neq q,
	\\
	\\
		\binom{n}{k} p\,^{k(n-k)}, & q = p.
	\end{array}
	\right.
\end{equation}

Due to \eqref{eq:recurseq} we can show that \emph{$\Tseq$-nomial} coefficients satisfy binomial-like recurrence relation, i.e., for any $n,k\in\Nat$ we have
\begin{equation}\label{eq:recur}
	\fnomialT{n}{k} = p^{n-k} \fnomialT{n-1}{k-1} + q^k \fnomialT{n-1}{k}
\end{equation}
with initial values $\fnomialT{n}{0} = \fnomialT{n}{n} = 1$.

% % % % % % % % % % % % % % % % % % % % % % % % % % % % % % % % % % 
\section{Generating functions of $\Tseq$-nomial coefficients \label{sec:gf}}

In the sequel, $\Tseq$ stand for a \emph{tileable} sequence $\TseqPQ{p}{q}$.

\begin{theoremn} \label{th:gf}
Let $\mathcal{A}_n(x)$ and $\mathcal{B}_n(x)$ be ordinary generating functions
\begin{subequations}
\begin{eqnarray}
	\label{eq:genA}
	\mathcal{A}_n(x) = \sum_{k\geq 0} (-1)^k q^{\binom{k}{2}} p^{\binom{k}{2}}  \fnomialT{n}{k} x^k,
	\\
	\label{eq:genB}
	\mathcal{B}_n(x) = \sum_{k\geq 0} \fnomialT{n+k-1}{k} x^k.
\end{eqnarray}
\end{subequations}
Then $\mathcal{A}_n(x)$ and $\mathcal{B}_n(x)$ are given as follows
\begin{subequations}
\begin{eqnarray}
	\label{eq:genAexpl}
	\mathcal{A}_n(x) = \prod_{i=1}^n \left( 1 - q^{(i-1)} p^{(n-i)}\, x \right)
	\\
	\label{eq:genBexpl}
	\mathcal{B}_n(x) = \prod_{i=1}^n \frac{1}{\left( 1 - q^{(i-1)} p^{(n-i)}\, x \right)}
\end{eqnarray}
\end{subequations}
with $\mathcal{A}_0(x) = 1$ and $\mathcal{B}_0(x) = 0$.
\end{theoremn}

\begin{Proof}
It is a simple exercise using ordinary generating functions. We apply recurrence \eqref{eq:recur} to \eqref{eq:genA} and \eqref{eq:genB} to see that
$$
\begin{array}{ll}
	\mathcal{A}_n(x) =  \left( 1 - p^{(n-1)} x \right) \mathcal{A}_{n-1} (q\,x), \quad &  \mathcal{A}_0(x) = 1; \\
	\\
	\mathcal{B}_n(x) =  \frac{1}{\left(1-p^{(n-1)} x\right)} \mathcal{B}_{n-1} (q\,x), & \mathcal{B}_0(x) = 0, \mathcal{B}_1(x) = \frac{1}{1-x}.
\end{array}
$$
And then it follows immediately that \eqref{eq:genAexpl} and \eqref{eq:genBexpl} hold by induction.
\end{Proof}

We can  infer also another form of the generating function $\mathcal{A}_n(x)$ of \emph{$\Tseq$-nomial} coefficients. Let $\mathcal{C}_n(x)$ be a generating function defined as follows
\begin{equation}
	\mathcal{C}_n(x) = \sum_{k\geq 0} (-1)^k q^{\binom{k}{2}} p^{\binom{n-k}{2}} \fnomialT{n}{k}.
\end{equation}
Then simple calculation using recurrence \eqref{eq:recur} yields
\begin{equation}
	\nonumber
	\mathcal{C}_n(x) = p^{n-1}(1-x) \mathcal{C}_{n-1}(x p/q)
\end{equation}
which immediately results in
\begin{equation}
	\mathcal{C}_n(x) = \prod_{i=1}^n \left( p^{(i-1)} - q^{(i-1)}x \right)
\end{equation}
for $n\geq 1$ and $\mathcal{C}_0(x) = 1$.

\begin{corol} \label{col:lambda}
Let $k,n \in\Nat$ be given. Then the value of $\Tseq$-nomial coefficient is equal
\begin{equation}
	\label{eq:otherform}
	\fnomialT{n+k-1}{k} = \sum_{1\leq b_1\leq \cdots \leq b_k\leq n}
	\lambda_{b_1}^n \lambda_{b_2}^n \cdots \lambda_{b_k}^n,
\end{equation}
where $\lambda_{i}^n = q^{(i-1)} p^{(n-i)}$ for $i=1,2,\ldots,n$.
\end{corol}

\begin{Proof}
	Let us consider generating function \eqref{eq:genBexpl}. From Theorem \ref{th:gf} we have
\begin{align}
	\nonumber
	\fnomialT{n+k-1}{k} &= [x^k] \prod_{i=1}^n \frac{1}{\left( 1 - q^{(i-1)} p^{(n-i)}\, x \right)} 
	= [x^k] \prod_{i=1}^n \sum_{j\geq 0} q^{(i-1)j} p^{(n-i)j} x^j
	\\
	\nonumber
	&= [x^k] \sum_{k\geq 0} c_{n,k} x^k = c_{n,k}
\end{align}
where coefficients $c_{n,k}$ take a form
\begin{align}
	\nonumber
	c_{n,k} &= \sum_{\substack{a_1+a_2+\cdots+a_n = k \\ a_1,a_2,\ldots,a_n\geq 0}} 
	q^{(1-1)a_1} p^{(n-1)a_1} q^{(2-1)a_2} p^{(n-2)a_2} \cdots q^{(n-1)a_n} p^{(n-n)a_n}
	\\ 
	\nonumber
	&= \sum_{\substack{a_1+a_2+\cdots+a_n = k \\ a_1,a_2,\ldots,a_n\geq 0}} 
	\left( q^{0} p^{(n-1)} \right)^{a_1} \left(q^{1} p^{(n-2)}\right)^{a_2} \cdots \left( q^{(n-1)} p^{0} \right)^{a_n}.
\end{align}

\noindent Notice that there are at most $k$ non-zero variables $a_i$ in the sum $a_1+a_2+\cdots+a_n = k$, where $a_i\geq 0$, for $i=1,2,\ldots,n$. Therefore, for fixed $k$ we choose a multiset $B$ (with repetition allowed) of $k$ these variables $a_i$ that are non-zero, i.e. $B = \{b_1,b_2,\ldots,b_k\}$ where $1\leq b_i\leq n$ and $1\leq b_1\leq b_2\leq\cdots \leq b_k\leq n$. Hence

\begin{align}
	\nonumber
	c_{n,k} &= \sum_{1\leq b_1\leq b_2\leq \cdots \leq b_k\leq n} 
	\left( q^{(b_1-1)} p^{(n-b_1)} \right) \left( q^{(b_2-1)} p^{(n-b_2)} \right) \cdots \left( q^{(b_k-1)} p^{(n-b_k)} \right)
	\\
	&= \sum_{1\leq b_1\leq b_2\leq \cdots \leq b_k\leq n} 
	\lambda_{b_1}^n \lambda_{b_2}^n \cdots \lambda_{b_k}^n
\end{align}
where $\lambda_i^n = q^{(i-1)} p^{(n-i)}$. Hence the thesis.
\end{Proof}

\begin{corol} \label{col:lambda2}
Let $k,n$ be natural numbers. Then the following hold
\begin{equation}
	\fnomialT{n}{k} q^{\binom{k}{2}} p^{\binom{k}{2}}
	= \sum_{1\leq b_1 < b_2 < \cdots < b_k \leq n}
		\lambda_{b_1}^n \lambda_{b_2}^n \cdots \lambda_{b_k}^n
\end{equation}
with $\lambda_i^n = q^{(i-1)} p^{(n-i)}$ for $i=1,2,\ldots,n$.
\end{corol}
\begin{Proof}
Consider generating function \eqref{eq:genAexpl}. It is easy to see that
\begin{equation}
	\nonumber
	\prod_{i=1}^n \left( 1 - \lambda_i^n x \right) 
	= \sum_{k\geq 0}^n  \left( \sum_{1\leq b_1 < b_2 < \cdots < b_k \leq n} \lambda_{b_1}^n \lambda_{b_2}^n \cdots \lambda_{b_k}^n  \right) (-1)^k x^k
\end{equation}
where $b_i$ denote indices of chosen factors $(1-\lambda_i^n x)$. Combining \eqref{eq:genA} with \eqref{eq:genAexpl} finishes the proof.
\end{Proof}

Since $p$ and $q$ are natural numbers, Corollary \ref{col:lambda} and Corollary \ref{col:lambda2} provide a combinatorial interpretation expressed in the language of object selections' from weighted boxes, where weight of $k$-th of $n$ boxes is given by $\lambda_k^n$. Let us sum up the above in the following corollary. Compare it also with the Konvalina \cite{konvalina,konvalina2} unified interpretation of Binomial, Gaussian coefficients and Stirling numbers.

\begin{corol} \label{col:interp}
Let $\Tseq$ be a tileable sequence $\TseqPQ{p}{q}$, where $p,q$ are natural numbers. Suppose we have $n$ boxes, such that $i$-th box contains $\lambda_i^n = q^{(i-1)} p^{(n-i)}$ distinguish balls, for $i=1,2,\ldots,n$.
Then the value of $\fnomialT{n+k-1}{k}$ is equal to the number of ways to select $k$ balls from $n$ boxes with box repetition allowed and $\fnomialT{n}{k} q^{\binom{k}{2}} p^{\binom{k}{2}}$ without repetition correspondingly.
\end{corol}

Note:  the appearance of nonnegative integers $p$ and $q$ gives another ways of formulation the combinatorial interpretation. For example, we can express it as a $k$-selection of pair of balls from two groups of boxes which size is dependent on $p$ and $q$ respectively etc.

\begin{corol}
The $\Tseq$-nomial coefficients designated by a sequence $\TseqPQ{p}{q}$ satisfy orthogonality relation, i.e.,
\begin{equation}
	\label{eq:orth1}
	\sum_{k=0}^s \fnomialT{n}{k} (-1)^{k} p^{\binom{k}{2}} q^{\binom{k}{2}} \fnomialT{n+s-k-1}{n-1} = 0
\end{equation}
or equivalent
\begin{equation}
	\label{eq:orth2}
	\sum_{k=0}^{n} \fnomialT{n+k-1}{k} (-1)^{n-k} p^{\binom{n-k}{2}} q^{\binom{n-k}{2}} \fnomialT{n}{k} = 0
\end{equation}
for any $n,s\in\Nat$.
\end{corol}

\begin{Proof}
Let $\mathcal{A}_n(x) = \sum_{k\geq 0} a_k x^k$ and $\mathcal{B}_n(x) = \sum_{k\geq 0} b_k x^k$ be generating functions of the form \eqref{eq:genAexpl} and \eqref{eq:genBexpl} accordingly. Observe that  $\mathcal{A}_n(x)\mathcal{B}_n(x) = 1$ for any $n\geq 1$. Therefore from the Cauchy product we have
$$
	\mathcal{A}_n(x) \cdot \mathcal{B}_n(x) = \sum_{s\geq 0} c_s x^s = 1
	\quad \Leftrightarrow \quad
	c_s = \sum_{k=0}^s a_k b_{s-k} = 0, \textrm{~for~} s\geq 1.
$$
The second form \eqref{eq:orth2} we can infer due to $\mathcal{B}_n(x)\mathcal{A}_n(x) = 1$.  Hence the thesis.
\end{Proof}

\begin{proposition}
Let $k,n$ be natural numbers. Suppose $p\neq q$, then the value of $\Tseq$-nomial coefficient is equal to
\begin{equation}
	\label{eq:explicit2}
	\fnomialT{n}{k} = \sum_{i=0}^{k} 
	(-1)^{(k-i)} \frac{(\mu_{i})^n}
	{\prod_{j=0}^{i-1} \left( \mu_{i} - \mu_{j} \right) \prod_{j=i+1}^{k} \left( \mu_j - \mu_i \right) }
\end{equation}
with $\mu_s = q^s p^{k-s}$ for $s = 0,1,2,\ldots,k$.
\end{proposition}

\begin{Proof}
We prove \eqref{eq:explicit2} using partial fractions expansion of the generating function \eqref{eq:genBexpl}, i.e.,

\begin{equation}
	\label{eq:proofexplicit1}
	\mathcal{B}_n(x) = \prod_{i=1}^n \frac{1}{(1 - q^{(i-1)} p^{(n-i)} x)} 
	=
	\sum_{i=1}^n \frac{a_i}{(1 - q^{(i-1)} p^{(n-i)} x)}
\end{equation}
If we multiply both sides of \eqref{eq:proofexplicit1} by $\prod_{j=1}^n (1 - q^{(j-1)} p^{(n-j)} x)$, then we can rewrite above as

$$
	1 = \sum_{i=1}^n a_i  \prod_{\substack{j=1 \\ j\neq i}}^n \left(1 - q^{(j-1)} p^{(n-j)} x\right).
$$
To find each $a_i$, let $x = 1/(q^{i-1} p^{n-i})$ for $i=1,2,\ldots,n$; and observe that all summands instead of $i$-th vanish, hence

$$
	a_i = \prod_{\substack{j=1 \\ j\neq i}}^n \frac{1}{ \left(1 - q^{(j-i)} p^{(i-j)}\right) }.
$$
Combining \eqref{eq:genB} with \eqref{eq:proofexplicit1} yields to

\begin{equation}
\nonumber
\begin{split}
	\fnomialT{n+k-1}{k} &
	= [x^k] \prod_{i=1}^{n} \frac{1}{\left( 1 - q^{(i-1)} p^{(n-i)} x \right)}
	= [x^k] \sum_{k\geq 0} \sum_{i=1}^n a_i \cdot q^{(i-1)k} p^{(n-i)k} x^k \\
	& = \sum_{i=0}^{n-1} \Bigg(  \prod_{\substack{j=0 \\ j\neq i}}^{n-1} \frac{1}{ \left(1 - q^{(j-i)} p^{(i-j)}\right) }  \Bigg) \left( q^{i} p^{n-i-1} \right)^{k}
\end{split}
\end{equation}
Let us make substitution $n = n'-k+1$ and use symmetry rule \eqref{eq:symmetry} of \emph{$\Tseq$-nomial} coefficients
\begin{equation}
\nonumber
\begin{split}
	\fnomialT{n+k-1}{k} & = \fnomialT{n'}{k} = \fnomialT{n'}{n'-k}
	= \sum_{i=0}^{n'\!-k} \Bigg(  \prod_{\substack{j=0 \\ j\neq i}}^{n'\!-k} \frac{1}{ \left(1 - q^{(j-i)} p^{(i-j)}\right) }  \Bigg) \left( q^{i} p^{n'-k-i} \right)^{k}.
\end{split}
\end{equation}
Replacing letter $k\geq 0$ by $(n'-k\,')\geq 0$ and next setting $\mu_i = q^i p^{(k'-i)}$ finishes the proof.
\begin{align}
	\nonumber
	\fnomial{n'}{k'} &= \sum_{i=0}^{k'} \Bigg(  \prod_{\substack{j=0 \\ j\neq i}}^{k'} \frac{1}{ \left(1 - q^{(j-i)} p^{(i-j)}\right) }  \Bigg) \left( q^{i} p^{k'-i} \right)^{n'-k'} 
	\\
	&= \sum_{i=0}^{k'} \frac{\left( q^{i} p^{k'-i} \right)^{n'}}
	{\prod_{\substack{j=0 \\ j\neq i}}^{k'} \left(q^i p^{(k'-i)} - q^{j} p^{(k'-j)}\right)}
	=
	\sum_{i=0}^{k'} \frac{(-1)^{(k'-i)} (\mu_i)^{n'}}
	{\prod_{j=0}^{i-1} \left(\mu_i - \mu_j\right) \prod_{j=i+1}^{k} \left(\mu_j - \mu_i\right)}
\end{align}
\end{Proof}

Note, that \eqref{eq:explicit2} is well defined for $n < k$, even for $n$ being a complex number.

\begin{exmp}
Take a sequence $\TseqPQ{q}{1}$ of Gaussian integers with arbitrary $q\neq 1$. Then $q$-binomial coefficient is given by
\begin{equation*}
	\fnomialF{n}{k}{q} = \sum_{i=0}^k (-1)^{i} 
	\frac { q^{(k-i)(n-i) - \binom{k-i}{2}}}
	{ \prod_{j=1}^i (q^j - 1) \prod_{j=1}^{k-i} (q^j - 1) }.
\end{equation*}
\end{exmp}

\begin{corol}
Let $k$ be a natural number. Then for any $p,q\in\mathbb{R}$, such that $p\neq q$ the following identity takes place
\begin{equation}
	\label{eq:equal1}
	\sum_{i=0}^{k} 
	(-1)^{(k-i)} \frac{(\mu_{i})^k}
	{\prod_{j=0}^{i-1} \left( \mu_{i} - \mu_{j} \right) \prod_{j=i+1}^{k} \left( \mu_j - \mu_i \right) } = 1
\end{equation}
where $\mu_i = q^i p^{(k-i)}$ for $i=0,1,2,\ldots,k$.
\end{corol}

\begin{Proof}
Indeed. Let $F$ be a sequence $\TseqPQ{p}{q}$. Then from definition of \emph{$F$-nomial} coefficients we have
$$
	\fnomial{k}{k} = 1
$$
On the other hand, setting $n=k$ in \eqref{eq:explicit2} completes the proof.
\end{Proof}

\begin{proposition}
Let $n, m, k\in\Nat$ and $\Tseq$ be a sequence $\TseqPQ{p}{q}$. Then the following relation holds
\begin{equation} \label{eq:conv}
	\fnomialT{n+m}{k} = \sum_{s=0}^k  p^{(m-k+s)}q^{(k-s)(n-s)} \fnomialT{n}{s} \fnomialT{m}{k-s}
\end{equation}
\end{proposition}

\begin{Proof}
It is an extension of the $q$-Vandermonde identity. Since \emph{$\Tseq$-nomial} satisfy \eqref{eq:otherform} we can prove \eqref{eq:conv} using Konvalina \cite{konvalina} way for $q$-binomial coefficients.
Assume $N = n+m-k+1$ and $\lambda_i^{N} = q^{(i-1)}p^{(N-i)}$ for $i=1,2,\ldots,N$. Then from \eqref{eq:otherform} we have that
\begin{equation}
	\nonumber
	\fnomialT{n+m}{k} = 
	\sum_{1\leq b_1\leq \cdots \leq b_k\leq N} 
	\lambda_{b_1}^N \lambda_{b_2}^N\cdots\lambda_{b_k}^N.
\end{equation}
We separate the above sum into $k+1$ sums, such that the first summation is over $\mathcal{S}um(0) := (n+1\leq b_1\leq\cdots\leq b_k\leq n+m-k+1)$, next $s$-th for $s=1,2,\ldots,k$ takes a form
\begin{equation}
	\nonumber
	\mathcal{S}um(s) := (1\leq b_1\leq \cdots \leq b_s \leq (n-s+1) \leq b_{s+1} \leq \cdots \leq b_k \leq n+m-k+1).
\end{equation}
In this convention we show that
\begin{equation}
	\nonumber
	\fnomialT{n+m}{k} = \sum_{s=0}^k 
	\sum_{\mathcal{S}um(s)} \lambda_{b_1}^N \lambda_{b_2}^N\cdots\lambda_{b_k}^N
	= 
	\sum_{s=0}^k 
	\sum_{\mathcal{S}um(s)} \left( q^{(b_1-1)} p^{(N-b_1)}\right) \cdots \left(q^{(b_k-1)} p^{(N-b_k)}\right).
\end{equation}
Consider now an $s$-th term of the above sum, and observe that $\mathcal{S}um(s)$ might be separated into two conditions
\begin{align}
	\nonumber
	\mathcal{S}um(s) &= 
	\big( 1\leq b_1\leq \cdots \leq b_s \leq n-s+1 \big) 
	\wedge 
	\big( n-s+1 \leq b_{s+1} \leq \cdots \leq b_k \leq n+m-k+1 \big)
\end{align}
Take $\mathcal{S}um(s) = \mathcal{S}um_1(s) \wedge \mathcal{S}um_2(s)$, where
\begin{equation}
	\nonumber
	\mathcal{S}um_1(s) := 
	\big( 1\leq b_1\leq \cdots \leq b_s \leq n-s+1 \big),
	\mathcal{S}um_2(s) :=
	\big( 1 \leq b'_{1} \leq \cdots \leq b'_{k-s} \leq m-(k-s)+1 \big).
\end{equation}
Therefore
\begin{equation}
	\nonumber
	\sum_{\mathcal{S}um(s)} \lambda_{b_1}^N \lambda_{b_2}^N \cdots \lambda_{b_k}^N
	= \sum_{\mathcal{S}um_1(s)} \lambda_{b_1}^N \cdots \lambda_{b_s}^N
	\sum_{\mathcal{S}um_2(s)} \lambda_{b'_1+n-s}^N \cdots \lambda_{b'_{k-s}+n-s}^N.
\end{equation}
Observe that for $i=1,2,\ldots,s$; the coefficient $\lambda_{b_i}^N$ is equal
\begin{equation}
	\nonumber
	\lambda_{b_i}^N = q^{(b_i-1)} p^{(n+m-k+1-b_i)} 
	= q^{(b_i-1)} p^{(n-s+1-b_i)} \cdot p^{(m+s-k)}
	= p^{(m+s-k)} \lambda_{b_i}^{n-s+1},
\end{equation}
while $\lambda_{b'_j+n-s}^N$ for $j=1,2,\ldots,k-s$ takes a form
\begin{equation}
	\nonumber
	\lambda_{b'_j+n-s}^N = q^{(b'_j-1)} p^{(m-(k-s)+1-b'_{j})} \cdot q^{(n-s)}
	= q^{(n-s)} \lambda_{b'_j}^{m-(k-s)+1}
\end{equation}
Hence due to equation \eqref{eq:otherform}
\begin{align}
	\nonumber
	\fnomialT{n+m}{k} &= \sum_{s=0}^k p^{(m+s-k)s} q^{(n-s)(k-s)} 
	\!\!\!\sum_{\mathcal{S}um_1(s)}\!\!\!\! \lambda_{b_1}^{n-s+1} \cdots \lambda_{b_s}^{n-s+1}
	\!\!\!\sum_{\mathcal{S}um_2(s)}\!\!\!\! \lambda_{b'_1}^{m-(k-s)+1} \cdots \lambda_{b'_{k-s}}^{m-(k-s)+1}
	\\
	&= \sum_{s=0}^k p^{(m+s-k)s} q^{(n-s)(k-s)} \fnomialT{n}{s} \fnomial{m}{k-s}.
\end{align}
\end{Proof}

% % % % % % % % % % % % % % % % % % % % % % % % % % % % % % % % % % 
\section{The Binomial-like theorem \label{sec:identity}}

In this section, we generalize generating function \eqref{eq:genAexpl} to the so called binomial-like theorem  (compare with Fib-binomial theorem in \cite{akk8}).

\begin{theoremn} \label{th:identity}
Let $n$ be a natural number $n\geq 1$. Then for any sequence $\TseqPQ{p}{q}$ we have
\begin{subequations} \label{eq:thiden}
\begin{align}
	\label{eq:thidenA}
	\prod_{i=1}^n (x+p^{(n-i)} q^{(i-1)} y)
	= \sum_{k\geq 0} \fnomialT{n}{k} q^{\binom{k}{2}} p^{\binom{k}{2}} y^{k}  x^{n-k} 
	\\
	\label{eq:thidenB}
	\prod_{i=0}^{n-1} (p^{i} x  +  q^{i} y)
	= \sum_{k\geq 0} \fnomialT{n}{k} q^{\binom{k}{2}} p^{\binom{n - k}{2}} y^{k} x^{n-k}
\end{align}
\end{subequations}
\end{theoremn}

\begin{Proof}
We first prove \eqref{eq:thidenA}. Let $\mathcal{C}_n(x,y)$ be an ordinary generating function
$$
	\mathcal{C}_n(x,y) = \sum_{k\geq 0} \fnomialT{n}{k} q^{\binom{k}{2}} p^{\binom{k}{2}} x^{n-k}y^k
$$
Applying recurrence relation \eqref{eq:recur} of \emph{$F$-nomial} coefficients for $k\geq 1$ yields 
\begin{equation}
	\nonumber
	\mathcal{C}_n(x,y) - x^n =
	\sum_{k\geq 1} \left\{ p^{n-k} \fnomialT{n-1}{k-1} + q^k \fnomialT{n-1}{k} \right\} q^{\binom{k}{2}} p^{\binom{k}{2}} x^{n-k} y^{k}.
\end{equation}
Let us observe that
$$
	\mathcal{C}_{n-1}(x,qy) = \sum_{k\geq 0} \fnomialT{n-1}{k} q^{\binom{k}{2}} p^{\binom{k}{2}} q^k y^k x^{n-1-k}.
$$
Now, it follows easily that $\mathcal{C}_0(x,y) = 1$ and for $n\geq 1$ function $\mathcal{C}_n(x,y)$ satisfies 
$$
	\mathcal{C}_n(x,y) = \left( x + p^{n-1}y \right) \mathcal{C}_{n-1}(x,qy).
$$
We can now proceed by induction over $n$ which completes the proof of \eqref{eq:thidenA}. In the same way we prove the second form \eqref{eq:thidenB}. However generating function $\mathcal{C}_{n}(x,y)$ is replaced by 
$$
	\mathcal{D}_n(x,y) = \sum_{k\geq 0} \fnomialT{n}{k} p^{\binom{n - k}{2}} x^{n-k}  q^{\binom{k}{2}} y^{k} 
$$
which follows to recurrence
$$
	\mathcal{D}_n(x,y) = \left( x + y \right) \mathcal{D}_{n-1}(p\,x,qy).
$$
with $\mathcal{D}_0(x,y) = 1$. And then it is easily seen to be \eqref{eq:thidenB}. Hence the thesis.
\end{Proof}

\noindent \textbf{Remark.} It is interesting to confront the Theorem \ref{th:identity} with logarithmic Fib-binomial theorem in \cite{akk8} valid for a very broad class of $F$-sequences - see Remark at the end of \cite{akk8}.

% % % % % % % % % % % % % % % % % % % % % % % % % % % % % % % % % % 
\section{Examples}

It is worth pointing out that the Kwa\'sniewski \emph{$F$-nomial} coefficients give us a unified interpretation of Binomial, Gaussian and Fibonomial coefficients. Recently it was showed \cite{md3} that these coefficients counts also bipartite, $k$-colored multi graphs and its inversion formula is related with the number of directed acyclic multi graphs. 

\vspace{0.2cm}
Throughout this part we shall consequently use the condition convention $n = k + m$.

% % % % % % % % % % % % % % % % % % % % % 
\subsection{Binomial coefficients}

If $F$ is a sequence of next natural numbers, or more general $n_F = 1_F\cdot n$ for any $1_F\in\Nat$. Then \emph{$F$-nomials} reduce to the basic Binomial coefficients. More precise, let $F$ be a sequence $\TseqPQ{p}{q}$ with $p=q=1$, then $n_F = k_F + m_F$ and
$$
	\fnomialF{n}{k}{\TseqPQ{1}{1}} \equiv \fnomialF{n}{k}{},
	\qquad
	\fnomialF{n}{k}{\TseqPQ{1}{1}} = \fnomialF{n-1}{k-1}{\TseqPQ{1}{1}} + \fnomialF{n-1}{k}{\TseqPQ{1}{1}}.
$$ 
Theorem \ref{th:identity} yields to well known binomial theorem
$$
	(x+y)^n = \sum_{k\geq 0} \fnomialF{n}{k}{\TseqPQ{1}{1}} x^{n-k}y^k
$$
From inversion formula \eqref{eq:inv} we obtain that
$$
	 \fnomialF{n}{k}{\TseqPQ{1}{1}}^{-1} =  (-1)^{n-k}\fnomialF{n}{k}{\TseqPQ{1}{1}}
$$
and then
$$
	(x - 1)^n = \sum_{k \geq 0} \fnomialF{n}{k}{\TseqPQ{1}{1}} (-1)^{n-k} x^k 
	\quad\Leftrightarrow\quad
	x^n = \sum_{k \geq 0} \fnomialF{n}{k}{\TseqPQ{1}{1}} (x - 1)^k.
$$

% % % % % % % % % % % % % % % % % % % % % 
\subsection{Gaussian coefficient}

If we take a sequence $\TseqPQ{p}{q}$ where $p=1$. Then we obtain well-known $q$-Calculus with Gaussian coefficients, i.e., let $n_F \equiv n_q = \frac{q^n - 1}{q - 1}$, then $n_F = k_F + q^k m_F$ and
\begin{equation}
	\fnomialF{n}{k}{\TseqPQ{1}{q}} \equiv \fnomialF{n}{k}{q} = \frac{n^{\underline{k}}_q}{k_q!},
	\qquad
	\fnomialF{n}{k}{\TseqPQ{1}{q}} = \fnomialF{n-1}{k-1}{\TseqPQ{1}{q}} + q^k \fnomialF{n-1}{k}{\TseqPQ{1}{q}}.
\end{equation}
with $q$-Binomial formula
\begin{equation}
	(1 - x)(1 - qx)\cdots(1-q^{n-1}x) = \sum_{k\geq 0} (-1)^k q^{\binom{k}{2}} \fnomialF{n}{k}{q} x^k
\end{equation}
and corresponding one
\begin{equation}
	\frac{1}{(1-x)(1-qx)\cdots(1-q^{k-1}x)} = \sum_{n\geq 0} \fnomialF{n+k}{k}{q} x^n.
\end{equation}
By Inversion formula \eqref{eq:inv} we immediately see that
$$
	\fnomialF{n}{k}{\TseqPQ{1}{q}}^{-1} = (-1)^{n-k}q^{\binom{n-k}{2}} \fnomialF{n}{k}{\TseqPQ{1}{q}}
$$
and consequently
$$
	\Phi_n(x) = \sum_{k \geq 0} \fnomialF{n}{k}{\TseqPQ{1}{q}} (-1)^{n-k} q^{\binom{n-k}{2}} x^k 
	\quad\Leftrightarrow\quad
	x^n = \sum_{k \geq 0} \fnomialF{n}{k}{\TseqPQ{1}{q}} \Phi_k(x)
$$
where $\Phi_n(x) = \prod_{s=0}^{n-1} (x - q^s)$.

% % % % % % % % % % % % % % % % % % % % % 
\subsection{Fibonomial coefficient}

In this section we show that the Fibonacci numbers define certain \emph{tileable} sequence $\TseqPQ{p}{q}$. Therefore in this case, the \emph{$\Tseq$-nomial} coefficients reduce to Fibonomial coefficients. 
Moreover an $m$-dimensional $F$-box $V_{k,n}$, such that 
$$
	V_{k,n} = [k_F]\times[(k+1)_F]\times \cdots \times [n_F] 
$$
where $s_F$ stays for $s$-th Fibonacci number might be filled with translates bricks $V_{1,m}$, where $m = n-k+1$ \cite{md1}.

\vspace{0.2cm}
Let $\alpha$ be natural number. If we set $\varphi = (\alpha \pm \sqrt{\alpha^2 + 4})/ 2$ then a sequence $F$ takes a form $n_F = \frac{1}{\sqrt{\alpha^2 + 4}}\left( \varphi_+^{\,n} - \varphi_-^{\,n} \right)$ with recurrence $n_F = \varphi_+^{\,m} k_F + \varphi_-^{\,k} m_F$, where $1_F = 1$ and $2_F = \alpha$. Note, that $\varphi_+, \varphi_-$ are real numbers. Therefore due to combinatorial interpretation we use another recurrence that Fibonacci numbers satisfy. Namely, let any $m,k$ such that $n= m+k$. Then
\begin{equation}
	n_F = (m-1)_F \cdot k_F + (k+1)_F \cdot m_F
\end{equation}
while $1_F = 1$ and $2_F = \alpha$. If we put $\alpha = 1$ we obtain sequence of Fibonacci numbers and Fibonomial coefficients
\begin{equation}
	\fnomialF{n}{k}{Fib} = (m-1)_F \fnomialF{n-1}{k-1}{Fib}  +  (k+1)_F \fnomialF{n-1}{k}{Fib}.
\end{equation}
Observe that $\varphi_+\cdot\varphi_- = -1$, then the generating function takes a form
\begin{equation}
	\prod_{s=1}^{n} \left( 1 - \varphi_-^{(s-1)} \varphi_+^{(n-s)} x\right) 
	= \sum_{k\geq 0} (-1)^{\binom{k+1}{2}} \fnomialF{n}{k}{Fib} x^k
\end{equation}
and
\begin{equation}
	\frac{x^k}{\prod_{s=0}^{k} \left( 1 - \varphi_+^{(k-s)} \varphi_-^{s}\, x \right)} =  \sum_{n\geq 0} \fnomialF{n}{k}{Fib} x^n.
\end{equation}

\subsection{Counting graphs $G(\alpha,n,k)$}

Labeled bipartite $\alpha$-multigraph $G(\alpha,n,k)$ is a bipartite graph with $n$ vertices ($k$ of them is in one of the two disjoint vertices' sets) with multiedges, such that any two vertices might be connected by at most $(\alpha-1)$ edges. We will denote by $\beta_{\alpha,n,k}$ the number of these graphs $G(\alpha,n,k)$.

\vspace{0.2cm} Let consider a sequence $\TseqPQ{p}{q}$ such that $p = q$. We can easily see that $n$-th term of $\TseqPQ{p}{p}$ is given $n_F = 1_F\cdot n p^{n-1}$ where $1_F\in\Nat$. For simplicity of notation we write $\TseqS{p}$ instead of $\TseqPQ{p}{p}$. It turns out that if $p\in\Nat$, then \emph{$F$-nomial} coefficients specified by sequence $\TseqS{p}$ define the number $\beta_{p,n,k}$, i.e.,

\begin{equation}
	\fnomialF{n}{k}{\TseqS{p}} = \fnomialF{n}{k}{} \cdot p^{k(n-k)} = \beta_{p,n,k}
\end{equation}

A recurrence relation of \emph{$F$-nomial} coefficients for a sequence $F = \TseqS{p}$ takes a form

$$
	\beta_{p,n,k} = \fnomialF{n}{k}{\TseqS{p}} = p^{n-k} \fnomialF{n-1}{k-1}{\TseqS{p}} + p^{k} \fnomialF{n-1}{k}{\TseqS{p}}
$$
with generating function

\begin{eqnarray}
	\left( 1 - p^{(n-1)}\, x \right)^n = \sum_{k\geq 0} (-1)^k p^{2\cdot\binom{k}{2}}  \fnomialF{n}{k}{\TseqS{p}} x^k,
	\\
	\frac{1}{\left( 1 - p^{(n-1)}\, x \right)^n} = \sum_{k\geq 0} \fnomialF{n+k-1}{k}{\TseqS{p}} x^k.
\end{eqnarray}

If $A_p(n)$ denote the number of labeled acyclic $p\,$-multi digraphs with $n$-nodes, i.e., such that any two nodes of it might be connected by at most $(p-1)$ arcs. Then following \cite{md3} we have

$$
	\fnomialF{n}{k}{\TseqS{p}}^{-1} = (-1)^{n-k} A_p(n-k) \fnomialF{n}{k}{\TseqS{p}}.
$$
while $A_p(0) = 1$.

\vspace{0.4cm}
\noindent \textbf{To be next...}

\vspace{0.2cm}
Corollary \ref{col:interp} provides a new combinatorial interpretation of \emph{$\Tseq$-nomial} coefficients. It is also very similar to the Konvaline unified interpretation of Binomial, Gaussian coefficients and Stirling numbers of two kinds. However the recurrence relation of Konvalina numbers do not reduce to recurrence of \emph{$\Tseq$-nomial} coefficients.
We expect that there is a way to generalize both of them.

\vspace{0.4cm}
\noindent \textbf{Acknowledgments}

\noindent I would like to thank my Professor A. Krzysztof Kwa\'sniewski for his very helpful comments, suggestions and improvements of this note.

\end{document}